\documentclass[a4paper,10pt]{article}
\usepackage{amsmath,amsthm,latexsym,makeidx,amscd,amssymb}

\numberwithin{equation}{subsection}

\newtheorem{prop}{Proposition}[section]
\newtheorem{lem}[prop]{Lemma}
\newtheorem{ddd}[prop]{Definition}
\newtheorem{theorem}[prop]{Theorem}

\newcommand{\ind}{\mathop{\mbox{\rm ind}}}

\newcommand{\dom}{\mathop{\rm dom}}
\newcommand{\Dom}{\mathop{\rm Dom}}
\newcommand{\Se}{\mathfrak S}
\newcommand{\SF}{\mathfrak{SF}}
\newcommand{\F}{\mathfrak{F}}
\newcommand{\ev}{{\rm ev}}

\newcommand{\Tr}{\mathop{\rm Tr}}

\newcommand{\C}{C^{\infty}}

\newcommand{\impl}{\Rightarrow}

\newcommand{\ve}{\varepsilon}

\newcommand{\A}{\mathcal A}
\newcommand{\ra}{\partial}

\DeclareMathOperator{\supp}{supp}
\DeclareMathOperator{\spfl}{sf}
\DeclareMathOperator{\Ran}{Ran}

\def\bbbr{{\rm I\!R}} 
\def\bbbn{{\rm I\!N}} 

\def\bbbc{{\rm I\!C}}

\def\bbbq{{\mathchoice {\setbox0=\hbox{$\displaystyle\rm Q$}\hbox{\raise
0.15\ht0\hbox to0pt{\kern0.4\wd0\vrule height0.8\ht0\hss}\box0}}
{\setbox0=\hbox{$\textstyle\rm Q$}\hbox{\raise
0.15\ht0\hbox to0pt{\kern0.4\wd0\vrule height0.8\ht0\hss}\box0}}
{\setbox0=\hbox{$\scriptstyle\rm Q$}\hbox{\raise
0.15\ht0\hbox to0pt{\kern0.4\wd0\vrule height0.7\ht0\hss}\box0}}
{\setbox0=\hbox{$\scriptscriptstyle\rm Q$}\hbox{\raise
0.15\ht0\hbox to0pt{\kern0.4\wd0\vrule height0.7\ht0\hss}\box0}}}}

\def\bbbz{{\mathchoice {\hbox{$\sf\textstyle Z\kern-0.4em Z$}}
{\hbox{$\sf\textstyle Z\kern-0.4em Z$}}
{\hbox{$\sf\scriptstyle Z\kern-0.3em Z$}}
{\hbox{$\sf\scriptscriptstyle Z\kern-0.2em Z$}}}}

\def\bbbc{{\mathchoice {\setbox0=\hbox{$\displaystyle\rm C$}\hbox{\hbox
to0pt{\kern0.4\wd0\vrule height0.9\ht0\hss}\box0}}
{\setbox0=\hbox{$\textstyle\rm C$}\hbox{\hbox
to0pt{\kern0.4\wd0\vrule height0.9\ht0\hss}\box0}}
{\setbox0=\hbox{$\scriptstyle\rm C$}\hbox{\hbox
to0pt{\kern0.4\wd0\vrule height0.9\ht0\hss}\box0}}
{\setbox0=\hbox{$\scriptscriptstyle\rm C$}\hbox{\hbox
to0pt{\kern0.4\wd0\vrule height0.9\ht0\hss}\box0}}}}

\title{A new topology on the space of unbounded selfadjoint operators and the spectral flow}
\author{Charlotte Wahl\footnote{This research was supported by a grant of AdvanceVT.}}
\date{}

\begin{document}
\maketitle

\begin{abstract}
We define a new topology, weaker than the gap topology, on the space of selfadjoint unbounded operators on a separable Hilbert space. We show that the subspace of selfadjoint Fredholm operators represents the functor $K^1$ from the category of compact spaces to the category of abelian groups and prove a similar result for $K^0$. We define the spectral flow of a continuous path of selfadjoint Fredholm operators generalizing the approach of Booss-Bavnek--Lesch--Phillips. 
\end{abstract}

\section*{Introduction}

The space of bounded Fredholm operators on a separable Hilbert space endowed with the norm topology is a classifying space for the functor $K^0$ from the category of compact spaces to the category of abelian groups \cite{ja}\cite{a}. The index map realizes an isomorphism between the $K$-theory of a point and $\bbbz$. Furthermore a particular connected component of the space of selfadjoint bounded Fredholm operators with the norm topology represents the functor $K^1$ \cite{as}. An isomorphism $K^1(S^1) \cong \bbbz$ is given by the spectral flow, which was introduced in \cite{aps}. 

These results can be applied to unbounded Fredholm operators by using the bounded transform $D \mapsto D(1+D^*D)^{-\frac 12}$. However, since many important geometric applications involve unbounded operators, it is more convenient to work directly with the space of unbounded selfadjoint Fredholm operators. The gap topology on the space of unbounded selfadjoint operators is the weakest topology such that the maps $D \mapsto (D \pm i)^{-1}$ are continuous. Gap continuity is weaker than continuity of the bounded transform. Booss-Bavnek--Lesch--Phillips defined the spectral flow for gap continuous paths \cite{blp} and Joachim proved that the space of unbounded selfadjoint Fredholm operators endowed with the gap topology is a classifying space for $K^1$ and the space of Fredholm operators with the subspace topology (see \S 1) is a classifying space for $K^0$ \cite{j}. 
 
In the first part of this paper we define and study a new topology on the space of unbounded selfadjoint operators. In this topology a path $(D_t)_{t \in [0,1]}$ is continuous if and only if the resolvents $(D_t \pm i)^{-1}$ depend in a strongly continuous way on $t$ and if there is an even function $\phi \in \C_c(\bbbr)$ with $\supp \phi=[-\ve,\ve]$ and $\phi'|_{(-\ve,0)}>0$ for some $\ve>0$ such that $\phi(D_t)$ is continuous in $t$. This topology is weaker than the gap topology. Compared with the latter it has some additional useful properties: The bounded transform of a continuous path is again continuous. If $(D_t)_{t \in [0,1]}$ is a continuous path of Fredholm operators and $(U_t)_{t \in [0,1]}$ is a strongly continuous path of unitaries, then $(U_tD_tU_t^*)_{t \in [0,1]}$ is again a continuous path of Fredholm operators. We show that the space of selfadjoint Fredholm operators endowed with this topology represents $K^1$ and the space of Fredholm operators with the subspace topology represents $K^0$. Furthermore we illustrate with an example that  families of Fredholm operators that are continuous with respect to this topology but not gap continuous arise naturally from differential operators on noncompact manifolds. Along the way we indicate how these results generalize to regular Fredholm operators on a Hilbert $C^*$-module. 

In the second part we define and study the spectral flow of a continuous path of selfadjoint Fredholm operators generalizing the approach of Booss-Bavnek--Lesch--Phillips and relate it to the winding number. The definition of the spectral flow given here is for paths with invertible endpoints equivalent to the definition of the noncommutative spectral flow in \cite{wa} applied to a separable Hilbert space. However, in \cite{wa} we used the theory of Hilbert $C^*$-modules in an essential way. One aim of this paper is to recover the results of \cite{wa} for a Hilbert space using classical functional analysis. We refer to \cite{wa} for applications.

\section{A new topology on the space of unbounded selfadjoint operators}

\label{class}

Let $H$ be a separable Hilbert space.

Recall that a closed densely defined operator $D$ on $H$ is called Fredholm if its bounded transform $F_D:=D(1+D^*D)^{-\frac 12}$ is Fredholm.

We denote the set of selfadjoint unbounded operators on $H$ by $S(H)$ and the set of selfadjoint unbounded Fredholm operators on $H$ by $SF(H)$.

As usual, $B(H)$ is the space of bounded operators on $H$ endowed with the norm topology and $K(H)$ is the subspace of compact operators.

Throughout let $B$ be a compact space. 

For a Banach space $V$ we denote by $C(B,V)$ the Banach space of continuous functions from $B$ to $V$ equipped with the supremum norm. We write $C(\bbbr)$ for $C(\bbbr,\bbbc)$. For $b \in B$ the evaluation map is $\ev_{b}:C(B,V) \to V,~f \mapsto f(b)$.

For a map $D:B \to S(H)$ we define 
$$\Dom D:=\{f \in C(B,H)~|~ f(b) \in \dom D(b) \mbox{ for all } b \in B \mbox{ and } Df \in C(B,H)\} \ .$$
Here $Df:B \to H$ is defined as $b \mapsto D(b)f(b)$. 

First we note some useful facts about the functional calculus of selfadjoint operators.

\begin{prop} 
\label{regul}
Let $D:B \to S(H)$ be a map.

The following conditions are equivalent:
\begin{enumerate}
\item At each point $b \in B$ the set $\ev_{b}(\Dom D) \subset \dom D(b)$ is a core for $D(b)$.
\item The resolvents $(D(b) \pm i)^{-1}$ depend in a strongly continuous way on $b \in B$.
\item For each $\phi \in C(\bbbr)$ the operator $\phi(D(b))$ depends in a strongly continuous way on $b \in B$.
\end{enumerate}
\end{prop}

\begin{proof}
Set $R_{\lambda}(b)=(D(b)+ \lambda)^{-1}$.

We show that (1) implies (2): Let $\lambda=\pm i$. Since $R_{\lambda}(b)$ is uniformly bounded, it is enough to prove that $\Dom R_{\lambda}$ is dense in $C(B,H)$. Let $f\in C(B,H)$ and let $\ve>0$. The assumption implies that the set $\ev_b((D+\lambda)(\Dom D))$ is dense in $H$ for any $b \in B$. Hence by compactness there is a finite open covering $\{U_j\}_{j \in I}$ of $B$ and functions $g_j \in \Dom D,~j \in I,$ such that $\|(D(b) + \lambda)g(b)_j-f(b)\| < \ve$ for all $b \in U_j$. Let $\{\chi_j\}_{j \in I}$ be a partition of unity subordinate to the covering $\{U_j\}_{j \in I}$ and set $f_j(b)=(D(b) + \lambda)g_j(b)$. Then $\sum_{j \in I}\chi_jf_j \in \Dom R_{\lambda}$ and $\|f-\sum_{j \in I}\chi_jf_j\| < \ve$. 

(2) $\impl$ (3):  Let $\phi \in C(\bbbr)$. Since the algebra generated by the functions $(x + i)^{-1}$ and $(x-i)^{-1}$ is dense in $C_0(\bbbr)$, the assertion holds for all $\psi \in C_0(\bbbr)$, in particular for $\psi(x)=\phi(x)(x+i)^{-1}$. Hence $\phi(D)f \in C(B,H)$ for $f \in R_iC(B,H)$. By a similar argument as above, (2) implies that $R_i C(B,H)$ is dense in $C(B,H)$. Since $\phi(D(b))$ is uniformly bounded, this implies the assertion.

(3) $\impl$ (2) $\impl$ (1) is clear. 
\end{proof}

\begin{lem}
\label{compcon}
Let $D:B \to S(H)$ be a map such that the resolvents $(D(b) \pm i)^{-1}$ depend in a strongly continuous way on $b \in B$. Assume that for each $b \in B$ there is given a symmetric operator $K(b)$ with $\dom D(b) \subset \dom K(b)$ such that $K(b)(D(b) + i)^{-1}$ is compact and depends continuously on $b$. Then 
for each $\phi \in C_0(\bbbr)$ $$\phi(D)- \phi(D+K) \in C(B,K(H)) \ .$$
\end{lem}

\begin{proof}
It is enough to prove the assertion for the functions $(x\pm i)^{-1}$. Since $\Dom D=\Dom (D+K)$, the previous proposition implies that $(D(b)+K(b) \pm i)^{-1}$ depends in a strongly continuous way on $b$.

Hence $(D(b) +K(b) \pm i)^{-1}-(D(b) \pm i)^{-1}=-(D(b) +K(b) \pm i)^{-1}K(b)(D(b) \pm i)^{-1}$ is compact and depends continuously on $b \in B$.
\end{proof}

\begin{lem}
\label{rescont}
Let $X$ be a topological space.
Let $D:X \to S(H)$ be a map such that the resolvents $(D(x) \pm i)^{-1} \in B(H)$ depend continuously on $x \in X$. Then $\phi(D(x))$ depends continuously on $x$ for any $\phi \in C_0(\bbbr)$. 
\end{lem}

\begin{proof}
This follows again from the fact that the functions $(x + i)^{-1}$ and $(x-i)^{-1}$ generate a dense subalgebra of $C_0(\bbbr)$.
\end{proof}

In particular, if $D:X \to B(H)$ is a continuous map such that $D(x)$ is selfadjoint for each $x \in X$, then $f(D(x))$ depends continuously on $x$ for all $f \in C(\bbbr)$.

Recall that the gap topology on $S(H)$ is the weakest topology such that the maps $$S(H) \to B(H),~D \mapsto (D+i)^{-1} \ ,$$
$$S(H) \to B(H),~ D \mapsto (D-i)^{-1}$$ 
are continuous. We denote by $S(H)_{gap}$ resp. $SF(H)_{gap}$ the set $S(H)$ resp. $SF(H)$ equipped with the gap topology. We refer to \cite{blp} for its properties.

In the following we introduce a new topology on $S(H)$.

Let $\phi \in \C_c(\bbbr)$ be an even function with $\supp \phi = [-1,1]$ and with $\phi'(x)>0$ for $x \in (-1,0)$. Define $\phi_n \in \C_c(\bbbr)$ by $\phi_n(x):=\phi(nx)$ for $n \in \bbbn$. 

Let $\Se_n(H)$ be the set $S(H)$ endowed with the weakest topology such that the maps 
$$\Se_n(H) \to H,~D \mapsto (D+i)^{-1}x \ ,$$
$$\Se_n(H) \to H,~D \mapsto (D-i)^{-1}x \ ,$$
$$\Se_n(H) \to B(H),~ D \mapsto \phi_n(D)$$
are continuous for all $x \in H$.

For any even function $\psi \in C_c(\bbbr)$ with $\supp \psi \in (-\frac 1n, \frac 1n)$ there is $g \in C_c(\bbbr)$ with $g(0)=0$ such that $\psi=g \circ \phi_n$. Hence $\Se_n(H) \to B(H),~D \mapsto \psi(D)$ is continuous. We will often make use of this property. It implies that the inclusion $\Se_m(H) \to \Se_n(H)$ is continuous for $m \le n$. Let $\Se(H)$ be the set $S(H)$ endowed with the direct limit topology.

Define $$SF_n(H):=\{D \in SF(H)~|~ \phi_n(D) \in K(H) \}$$ and denote by $\SF_n(H)$ the set $SF_n(H)$ endowed with the subspace topology of $\Se_n(H)$.
Let $\SF(H)$ be the inductive limit of the spaces $\SF_n(H)$. 
An operator $D \in S(H)$ is Fredholm if and only if $F_D$ is invertible in $B(H)/K(H)$, and this is equivalent to $\phi_n(D) \in K(H)$ for $n$ big enough. 
Hence the underlying set of $\SF(H)$ is $SF(H)$.

If $D:B \to \Se(H)$ is continuous, then $f(D):B \to \Se(H)$ is continuous for any odd non-decreasing continuous function $f:\bbbr \to \bbbr$ with $f^{-1}(0)=\{0\}$. This can be seen as follows: Assume that $D: B \to \Se_n(H)$ is continuous. Since $(f \pm i)^{-1} \in C(\bbbr)$, we get from Prop. \ref{regul} that $(f(D) \pm i)^{-1}x:B \to H$ is continuous for any $x \in H$. Furthermore for $m$ big enough $\supp(\phi_m\circ f) \subset (-\frac 1n,\frac 1n)$, hence  $\phi_m(f(D)):B \to B(H)$ is continuous.

In particular the bounded transform $B \to \Se(H),~b \mapsto F_{D(b)}$ is continuous. The example of Fuglede presented in \cite{blp} shows that the bounded transform of a gap continuous family is in general not gap continuous. 

We need the following technical lemmata.

\begin{lem}
\label{compcont}
Assume that $D:B\to \SF_n(H)$ is continuous. Then for $\psi \in C_c(\bbbr)$ with $\supp \psi \subset (-\frac 1n,\frac 1n)$ we have that $\psi(D):B \to K(H)$ is continuous.
\end{lem}

\begin{proof} This follows from an elementary argument in the theory of Hilbert $C^*$-modules:

Let $B(C(B,H))$ be the algebra of strongly continuous families of bounded operators on $H$ with parameter space $B$ and with adjoint depending in a strongly continuous way on the parameter. Endowed with the supremum norm this is a $C^*$-algebra and $C(B,K(H))$ defines a closed ideal in $B(C(B,H))$.
 
Let $\pi:B(C(B,H)) \to B(C(B,H))/C(B,K(H))$ be the projection. Let $g \in C(\bbbr)$ with $g(0)=0$ be such that $\psi^2=g \circ \phi_n$. We have that $\pi(\phi_n(D))=0$, hence $\pi(\psi(D))^2=g(\pi(\phi_n(D)))=0$. Since $\pi(\psi(D))$ is selfadjoint in the $C^*$-algebra $B(C(B,H))/C(B,K(H))$, it follows that $\pi(\psi(D))=0$, hence $\psi(D) \in C(B,K(H))$. 
\end{proof} 

\begin{lem}
\label{compcont2}
If $(F_b)_{b \in B}$ is a strongly continuous family of bounded selfadjoint operators such that $(b\mapsto F_b^2-1) \in C(B,K(H))$, then for any function $\phi \in C(\bbbr)$ with $\phi(1)=\phi(-1)=1$ we have that $(b \mapsto \phi(F_b)-1) \in C(B,K(H))$.
\end{lem}

\begin{proof} The argument is similar to the proof of Lemma \ref{compcont}. We use its notation.

We have that $\pi((F_b)_{b \in B})^2=1$; hence the spectrum of $\pi((F_b)_{b \in B})$ is a subset of $\{-1,1\}$; thus $\phi(\pi((F_b)_{b \in B}))-1=0$. Since $\phi(\pi((F_b)_{b \in B}))=\pi((\phi(F_b))_{b \in B})$, it follows that $(b \mapsto \phi(F_b)-1) \in C(B,K(H))$.
\end{proof}

The following two properties of the space $\SF(H)$ are useful:

\begin{itemize}
\item Assume that $D:B \to \SF(H)$ is continuous and that $B \ni b \mapsto U(b)$ is a map with values in the group of unitaries of $B(H)$ such that $U(b)$ depends in a strongly continuous way on $b$. Then $UDU^*: B \to \SF(H)$ is continuous. 
\item If $D:B \to \SF(H)$ is continuous, then $f(D):B \to \SF(H)$ is continuous for any non-decreasing continuous function $f: \bbbr \to \bbbr$ with $f^{-1}(0)=\{0\}$. 
\end{itemize}

The first property follows from the fact that the composition of a continuous family of compact operators with of a strongly continuous family of bounded operators is again continuous if the parameter space is compact. Furthermore since $U$ is bounded below, the adjoint depends also in a strongly continuous way on $b$.

Note that the second property does not assume the function to be odd. Taking the Lemma \ref{compcont} into account one proves the property analogously to the corresponding one for $\Se(H)$ from above.

\begin{lem}
Let $D:B \to \SF(H)$ be continuous. Then there is an odd non-decreasing function $\chi \in C(\bbbr)$ with $\chi^{-1}(0)=\{0\}$ and $\lim_{x \to \infty}\chi(x)=1$ such that $\chi(D)^2-1:B \to K(H)$ is continuous.
\end{lem}

\begin{proof}
There is $n \in \bbbn$ such that
 $D:B \to \SF_n(H)$ is continuous. Then any non-decreasing $\chi\in C(\bbbr)$ with $\chi^{-1}(0)=\{0\}$ and such that $\supp (\chi^2-1) \subset (-\frac 1n,\frac 1n)$ works.
\end{proof}

\begin{ddd}
Let $D:B \to \SF(H)$ be continuous. Then a function $\chi$ fulfilling the conditions of the previous lemma is called a {\rm normalizing function} for $D$.
\end{ddd}

The name ``normalizing function'' is borrowed from \cite{hr}. The definition in \cite{hr} is different since it applies to a different class of operators, but the underlying idea is the same.

The definition of the spaces $\Se(H)$ and $\SF(H)$ generalizes in a straightforward way to the case where $H$ is a Hilbert $C^*$-module. In this case we assume the unbounded operators to be regular.

The spaces $\Se(H)$ of $S(H)_{gap}$ share many properties as we will see in the following. We omit some details since the arguments resemble those in \cite{blp}. 

First we note that $\SF(H)$ is path-connected since $SF(H)_{gap}$ is path-connected by \cite[Th. 1.10]{blp}.

Let $D_0 \in S(H)$. For $n\in \bbbn$ and $\ve>0$ we define $$U(\ve,n,D_0):=\{D \in S(H)~|~\|\phi_n(D)-\phi_n(D_0)\| < \ve\} \ .$$
This is an open neighbourhood of $D_0$ in $\Se_n(H)$. 

Let $(a,b)\subset \bbbr$ be in the resolvent set of $\phi_n(D_0)$. Then there is $\ve>0$ such that $(a,b)$ is in the resolvent set of $\phi_n(D)$ for all $D \in U(\ve,n,D_0)$. Hence $\phi_n^{-1}((a,b))$ is in the resolvent set of $D$ for all $D \in U(\ve,n,D_0)$. Furthermore if $\mu \in \phi_n^{-1}((a,b)),~\mu>0$, then also $-\mu \in \phi_n^{-1}((a,b))$ and $U(\ve,n,D_0) \to B(H),~D \mapsto 1_{[-\mu,\mu]}(D)$ is continuous.

This implies the following lemma, which will be used for the definition of the spectral flow:

\begin{lem}
\label{opensp}
If $D_0 \in SF_n(H)$ and $\mu \in (0,\frac 1n)$ is such that $\pm \mu$ is in the resolvent set of $D_0$, then there is $\ve>0$ such that $\pm \mu$ is in the resolvent set of $D$ for all $D \in U(\ve,n, D_0)$. Furthermore $1_{[-\mu,\mu]}(D)$ has finite-dimensional range for all $D \in U(\ve,n, D_0)$ and the map
$$\Se_n(H) \supset U(\ve,n, D_0) \to K(H),~ D \mapsto 1_{[-\mu,\mu]}(D)$$ is continuous.

In particular all operators in $U(\ve,n,D_0)$ are Fredholm.
\end{lem}

Note that for a given $D_0 \in SF_n(H)$ a $\mu$ fulfilling the assumption of the lemma always exists since the spectrum of $D_0$ near zero is discrete.

\begin{prop}
\begin{enumerate}
\item The identity induces a continuous map $S(H)_{gap} \to \Se(H)$.
\item The space $SF(H)$ is open in $\Se(H)$.
\item The identity induces a homeomorphism from $\Se(H) \cap SF(H)$ to $\SF(H)$.
\end{enumerate}
\end{prop}

\begin{proof}
The first assertion is a consequence of Lemma \ref{rescont}.

The second assertion follows from the previous lemma and the subsequent remark. Since the remark is in general wrong for a Hilbert $C^*$-module, we give another argument which also works for Hilbert $C^*$-modules: Let $D_0 \in SF_n(H)$ and let $\chi$ be a normalizing function for $D_0$ with $\supp (\chi^2-1) \subset (-\frac 1n,\frac 1n)$. Then $\chi(D_0)^2$ is invertible in $B(H)/K(H)$. Furthermore, since $\Se_n(H) \to B(H),~D \mapsto (\chi(D)^2-1)$ is continuous, also $\Se_n(H) \to B(H)/K(H),~D \mapsto  \chi(D)^2$ is continuous. Hence there is an open neighbourhood $U$ of $D_0$ in $\Se_n(H)$ such that $\chi(D)^2$ is invertible in $B(H)/K(H)$ for all $D \in U$. This implies that all $D \in U$ are Fredholm. 

The third assertion is clear.
\end{proof}

We denote the space of (not necessarily selfadjoint) Fredholm operators on $H$ by $F(H)$. We identify $F(H)$ with a subspace of $SF(H \oplus H)$ via the injection $$F(H) \to SF(H\oplus H),~D \mapsto \left(\begin{array}{cc} 0 & D^* \\ D & 0 \end{array}\right) \ .$$ 
Note that if $D \in F(H)$ and $f:\bbbr \to \bbbr$ is an odd non-decreasing continuous function with $f^{-1}(0)=\{0\}$, then $f(D) \in F(H)$ is well-defined.

The space $F(H)$ endowed with the subspace topology of $\SF(H\oplus H)$ is denoted by $\F(H)$. 

For topological spaces $X,Y$ we denote by $[X,Y]$ the set of homotopy classes of continuous maps from $X$ to $Y$.

\begin{theorem}
\label{classsp}
\begin{enumerate}
\item The space $\SF(H)$ represents the functor $B \mapsto K^1(B)$ from the category of compact spaces to the category of abelian groups.
\item The space $\F(H)$ represents the functor $B \mapsto K^0(B)$ from the category of compact spaces to the category of abelian groups. 
\end{enumerate}
\end{theorem}

\begin{proof}
We use the notation of \cite{j}: Let $KC_{sa}(H)$ (where $KC$ stands for ``Kasparov cycles'') be the space of selfadjoint bounded operators $F$ on $H$ with $\|F\|\le1$ and $F^2-1\in K(H)$ and endow it with the weakest topology such that the maps
$$KC_{sa}(H) \to H,~F \mapsto Fx \ ,$$ 
$$KC_{sa}(H) \to K(H),~ F \mapsto F^2-1$$ are continuous for all $x \in H$.
The inclusion $KC_{sa}(H) \to \SF(H)$ is continuous.

Let $KC(H)$ be the space of bounded operators $F$ such that $\|F\| \le 1$ and $F^*F-1,~FF^*-1\in K(H)$. Consider $KC(H)$ as a subspace of $KC_{sa}(H \oplus H)$ as above. 

By \cite[Theorem 3.4]{j}, which is based on results of Bunke--Joachim--Stolz, the space $KC_{sa}(H)$ represents the functor $K^1$ and the space $KC(H)$ represents $K^0$.

(1) Let $D:B \to \SF(H)$ be a continuous map. Let $\chi$ be a normalizing function for $D$ and let $\chi_t(x)=(1-t)x + t \chi(x)$. Then $B \to KC_{sa}(H),~b \mapsto \chi_1(D(b))$  and $[0,1]\times B \to \SF(H),~(t,b) \mapsto \chi_t(D(b))$ are continuous (here we use Prop. \ref{regul}). It follows that the map $[B,KC_{sa}(H)] \to [B,\SF(H)]$ induced by the inclusion $KC_{sa}(H) \to \SF(H)$ is surjective. 

For injectivity let $h:[0,1] \times B \to \SF(H)$ be a homotopy between continuous maps $B \to KC_{sa}(H),~b\mapsto h(i,b),~i=0,1$. Let $\chi$ be a normalizing function for $h$ such that $\chi(1)=1$ and let $\chi_t(x)=(1-t)x + t \chi(x)$. Since $\chi_t^2(1)-1= \chi_t^2(-1)-1=0$, Lemma \ref{compcont2} implies that the map $$B \to K(H),~b \mapsto \chi_t(h(i,b))^2-1$$ is continuous for $i=0,1$. Furthermore $\chi_t(h(i,b))^2-1$ is continuous in $t$ since $\chi_t^2-1$ depends continuously on $t$ in $C(\bbbr)$. Hence the map $$([0,1] \times \{0,1\}\times B) \cup (\{1\}\times [0,1] \times B) \to KC_{sa}(H),~(t,x,b) \mapsto \chi_t(h(x,b))$$ is continuous and defines a homotopy in $KC_{sa}(H)$ between $\chi_0(h(0,\cdot))=h(0,\cdot)$  and $\chi_0(h(1,\cdot))=h(1,\cdot)$.    
 
(2) The proof is analogous with the obvious modifications. 
\end{proof} 

It follows that $\pi_0(\F(H)) \cong \bbbz$. As usual an isomorphism is given by the index map.
The results in the following section will imply that an isomorphism $[S^1,\SF(H)] \to \bbbz$ is given by the spectral flow.

The proof of the previous proposition carries over to the case where $H$ is the standard Hilbert $\A$-module $H_{\A}$ of a unital $C^*$-algebra $\A$ implying that $\SF(H_{\A})$ is a representing space of the functor $B \mapsto K_1(C(B,\A))$ from the category of compact spaces to the category of abelian groups and $\F(H_{\A})$ is a representing space for $B \mapsto K_0(C(B,\A))$. The corresponding statements for $SF(H_{\A})_{gap}$ have been proven in \cite{j}. 

In the following we give two examples of maps with values in $SF(H)$ that are continuous in $\SF(H)$ but not gap continuous. Both arise from elliptic differential operators on a noncompact manifold. 

Let $H=L^2(\bbbr)$ and let $f \in \C(\bbbr)$ be nonconstant real-valued and bounded below by some $c>0$. Set $f_t(x)= f(tx)$ for $t \in [0,1]$. We define $D(t)$ on $L^2(\bbbr)$ to be the multiplication by $f_t$. The path $D:[0,1] \to SF(H)$ is not gap continuous at $t=0$, but it is continuous as a path in $\SF(H)$. 

Even if the resolvents are compact, they need not depend in a continuous way of $t$: Let $H=L^2(\bbbr,\bbbc^2)$.
Let $f \in \C_0(\bbbr)$ be a nonnegative function and let $g \in \C(\bbbr)$ with $g \ge 0$, $g(0)=1$, $g(1)=0$ and $g(x)=1$ for $|x| \ge 2$. Define $\psi_t(x):=\frac{g(tx)}{f(x)}+1$ for $t \in [0,1]$. Note that $\psi_t(x)$ is continuous in $t$ and $x$. Define $D(t)$ to be the closure of $$\left(\begin{array}{cc} 0 & \psi_t(1-\ra_x^2) \\ (1-\ra_x^2)\psi_t & 0 \end{array}\right):\C_c(\bbbr,\bbbc^2) \to L^2(\bbbr,\bbbc^2) \ .$$ Since $\frac{1}{\psi_t} \in C_0(\bbbr)$, the operator $\frac{1}{\psi_t}(1-\ra_x^2)^{-1}$ is compact on $L^2(\bbbr)$ for any $t$, hence $D(t)^{-1}$ is compact for any $t$. Furthermore $D(t)^{-1}$ is uniformly bounded. Thus $[0,1] \to \SF(H),~ t \mapsto D(t)$ is continuous. It is easy to check that $D(t)^{-1}$ is not continuous in $t$ at $t=0$. Hence $D$ is not gap continuous.

Note that these examples have in common that the coefficients are continuous as maps from $[0,1]$ to $C_{loc}(\bbbr)$ but not continuous (in the second example even not well-defined) as maps from $[0,1]$ to $C(\bbbr)$.

See \cite[\S 6]{wa} for criteria for the continuity  in $\SF(H)$ of families of elliptic operators on noncompact Riemannian manifolds and families of well-posed boundary value problems.

\section{Spectral flow}

In the following we generalize the definition of the spectral flow in \cite{blp}, which is based on the approach of \cite{p}, to continuous paths in $\SF(H)$. 

\begin{ddd}
Let $(D_t)_{t\in [a,b]}$ be a continuous path in $\SF(H)$ and assume that there is $\mu >0$ such that $\pm \mu$ is in the resolvent set of $D_t$ for all $t \in [a,b]$ and $1_{[-\mu,\mu]}(D_t)$ has finite-dimensional range. We define $$\spfl((D_t)_{t\in [a,b]})= \dim \Ran(1_{[0,\mu]}(D_b)) -\dim \Ran( 1_{[0,\mu]}(D_a)) \ .$$

If $(D_t)_{t\in [a,b]}$ is a general continuous path in $\SF(H)$, then we define its spectral flow by cutting the path into smaller pieces to which the previous situation applies and adding up the contributions of the pieces. (This is always possible by Lemma \ref{opensp} and the subsequent remark.)
\end{ddd}

Well-definedness can be proven as in \cite{p}. 

The spectral flow has the following properties:
\begin{enumerate}
\item It is additive with respect to concatenation of paths. 
\item For any non-decreasing continuous function $f:\bbbr \to \bbbr$ with $f^{-1}(0)=\{0\}$ 
$$\spfl((D_t)_{t\in [a,b]}) = \spfl((f(D_t))_{t \in [a,b]}) \ .$$
\item If $(U_t)_{t \in [a,b]}$ is a strongly continuous path of unitaries on $H$, then $$\spfl((U_tD_tU_t^*)_{t \in [a,b]})=\spfl((D_t)_{t \in [a,b]}) \ .$$
\item If $D_t$ is invertible for any $t \in [0,1]$, then $\spfl((D_t)_{t \in [a,b]})=0 \ .$
\item If $(D_{(s,t)})_{(s,t) \in [0,1]\times [a,b]}$ is a continuous family in $\SF(H)$ such that $D_{(s,a)}$ and $D_{(s,b)}$ are invertible for all $s \in [0,1]$, then
$$\spfl((D_{(0,t)})_{t\in [a,b]})=\spfl((D_{(1,t)})_{t \in [a,b]}) \ .$$ 
\item If $(D_{(s,t)})_{(s,t) \in [0,1]\times [a,b]}$ is a continuous family in $\SF(H)$ such that $D_{(s,a)}=D_{(s,b)}$ for all $s \in [0,1]$, then
$$\spfl((D_{(0,t)})_{t\in [a,b]})=\spfl((D_{(1,t)})_{t \in [a,b]}) \ .$$
\end{enumerate}

The proof of the first three properties is not difficult and is left to the reader.
The fourth property follows from the fact that by Lemma \ref{opensp} and by compactness of $[a,b]$ there is $\delta>0$ such that $[-\delta,\delta]$ is a subset of the resolvent set of $D_t$ for all $t \in [a,b]$ if $D_t$ is invertible for all $t \in [a,b]$.

The following proposition implies the last two properties, namely homotopy invariance:

\begin{prop}
If $(D_{(s,t)})_{(s,t) \in [0,1] \times [a,b]}$ is a continuous family in $\SF(H)$, then
$$\spfl((D_{(0,t)})_{t\in [a,b]}) +\spfl((D_{(s,b)})_{s \in [0,1]}) -\spfl((D_{(1,t)})_{t \in [a,b]})-\spfl((D_{(s,a)})_{s \in [0,1]})=0 \ .$$ 
\end{prop}

\begin{proof}
Let $n \in \bbbn$ be such that the family $(D_{(s,t)})_{(s,t) \in [0,1] \times [a,b]}$ is continuous in $\SF_n(H)$. 

If there is $\mu \in (0, \frac 1n)$ such that $\pm \mu$ is in the resolvent set of $D_{(s,t)}$ for all $(s,t) \in [0,1]\times [a,b]$, then  $1_{[-\mu,\mu]}(D_{(s,t)})$ has finite-dimensional range for all $(s,t)$ and the assertion follows from the definition of the spectral flow. 

In general we find, by compactness of $[0,1] \times [a,b]$ and by Lemma \ref{opensp}, an $n \in \bbbn$ such that each of the rectangles $[\tfrac{(m_1-1)}{n}, \tfrac{m_1}{n}] \times [a +(b-a)\tfrac{m_2-1}{n}, a+ (b-a) \tfrac{m_2}{n}]$ with $m_1,m_2=1,2 \dots n$ has the following property: There is a $\mu \in (0, \frac 1n)$ such that $\pm \mu$ is in the resolvent set of $D_{(s,t)}$ for all points $(s,t)$ of the rectangle. Hence for each of the rectangles an analogue of the formula holds by the previous argument. Since for fixed $n$ these rectangles constitute a subdivision of $[0,1]\times [a,b]$, the formula follows from the additivity of the spectral flow with respect to concatenation.
\end{proof}

We draw some conclusions in the following two propositions. See \cite[\S 3]{le} for similar results.

If $P,Q$ are projections with $P-Q \in K(H)$, then $QP:P(H) \to Q(H)$ is Fredholm with parametrix $PQ$. Let $$\ind(P,Q):= \ind (QP:P(H) \to Q(H)) \ .$$
It is well-known that $$\spfl((t(2P-1)+(1-t)(2Q-1))_{t \in [0,1]})=\ind(P,Q) \ .$$ 

\begin{prop}
\label{proj} 
Let $(P_t)_{t \in [0,1]},~(Q_t)_{t \in [0,1]}$ be strongly continuous paths of projections on $H$ such that $P_t-Q_t$ is compact and continuous in $t$. Then
$$\ind(P_0,Q_0)=\ind(P_1,Q_1) \ .$$
\end{prop}

\begin{proof}
First we prove that the family $(F_{(s,t)})_{(s,t) \in [0,1]^2}$ defined by $F_{(s,t)}:=t(2P_s-1)+(2t-1)(2Q_s-1)$ is continuous in $\SF(H)$: Clearly $F_{(s,t)}$ depends in a strongly continuous way on $(s,t)$. Hence, by Prop. \ref{regul}, the operators $(F_{(s,t)}\pm i)^{-1}$ depend in a strongly continuous way on $(s,t)$ as well. Furthermore $F_{(s,t)}-(2P_s-1)$ is a compact operator depending continuously on $(s,t)$. This and Lemma \ref{compcon} imply that $\phi_n(F_{(s,t)})-\phi_n((2P_s-1))$ is a compact operator depending continuously on $(s,t)$ for any $n \in \bbbn$. From $\phi_n((2P_s-1))=0$ it follows that $\phi_n(F_{(s,t)})$ is a compact operator depending continuously on $(s,t)$. This shows the continuity.

Now by homotopy
invariance $$\spfl((t(2P_0-1)+(1-t)(2Q_0-1))_{t \in [0,1]})=\spfl((t(2P_1-1)+(1-t)(2Q_1-1))_{t \in [0,1]}) \ .$$
\end{proof}

The following technical lemma, which is an immediate consequence of  \cite[Prop. 3.4]{le} and Lemma \ref{compcon}, is needed for the proof of the subsequent proposition:

\begin{lem}
Let $D \in S(H)$ and let $K$ be a symmetric operator with $\dom D \subset \dom K$ and such that $K(D+i)^{-1}$ is compact. Then $f(D+K)-f(D) \in K(H)$ for each function $f\in C(\bbbr)$ for which $\lim\limits_{x \to \infty}f(x)$ and  $\lim\limits_{x \to -\infty}f(x)$ exist.
\end{lem}

\begin{prop}
Let $(D_t)_{t \in [a,b]}$ be a continuous path in $\SF(H)$ with invertible endpoints and assume that there is a path of symmetric operators $(K_t)_{t \in [a,b]}$ with $\dom D_t \subset \dom K_t$ for all $t \in [a,b]$, such that $K_t(D_t + i)^{-1}$ is compact and continuous in $t$ and such that $(D_t+K_t)$ is invertible for each $t \in [a,b]$. 

Then \begin{eqnarray*}
\lefteqn{\spfl((D_t)_{t \in [a,b]})}\\
&=&\ind(1_{\ge 0}(D_b),1_{\ge 0}(D_b+K_b)) - \ind(1_{\ge 0}(D_a),1_{\ge 0}(D_a+K_a)) \ .
\end{eqnarray*}
\end{prop}

\begin{proof}
Let $n$ be such that $(D_t)_{t \in [a,b]}$ is a continuous path in $\SF_n(H)$.

Lemma \ref{compcon} implies that $\phi_n(D_t)-\phi_n(D_t+K_t)$ is compact and continuous in $t$. In particular $(D_t+K_t)_{t \in [a,b]}$ is a continuous path in $\SF_n(H)$. Since each $D_t+K_t$ is invertible, by property (4) of the spectral flow $$\spfl((D_t+K_t)_{t \in [a,b]})=0 \ .$$
  
Let $\psi \in C(\bbbr)$ with $\psi|_{(-\infty,\frac 13]}=0$ and $\psi|_{[\frac 23,\infty)}=1$.

By homotopy invariance and additivity with respect to concatenation the spectral flow of the path $(D_t)_{t \in [a,b]}$ equals the spectral flow of the path $({\tilde D}_t)_{t \in [a-1,b+1]}$ with ${\tilde D}_t=D_a + \psi(t-a+1)K_a$ for $t \in [a-1,a]$, ${\tilde D}_t=D_b+(1-\psi(t-b))K_b$ for $t \in [b,b+1]$ and ${\tilde D}_t=D_t+K_t$ for $t \in [a,b]$. Furthermore by additivity with respect to concatenation and since $\spfl((D_t+K_t)_{t \in [a,b]})=0$,
$$\spfl(({\tilde D}_t)_{t \in [a,b]})=\spfl(({\tilde D}_t)_{t \in [a-1,a]})+\spfl(({\tilde D}_t)_{t \in [b,b+1]}) \ .$$ 
We calculate $\spfl(({\tilde D}_t)_{t \in [a-1,a]})$: Let $\chi \in \C(\bbbr)$ be a normalizing function for $({\tilde D}_t)_{t \in [a-1,a]}$ such that $\chi(D_a)=2\cdot 1_{\ge 0}(D_a)-1$ and $\chi(D_a+K_a)=2\cdot 1_{\ge 0}(D_a+K_a)-1$. By the previous lemma $\chi(D_a)-\chi(D_a+K_a) \in K(H)$. Then 
\begin{eqnarray*}
\spfl(({\tilde D}_t)_{t \in [a-1,a]})&=& \spfl((\chi({\tilde D}_t))_{t \in [a-1,a]})\\
&=&\spfl(((1-t)\chi(D_a)+t\chi(D_a+K_a))_{t \in [0,1]})\\
&=&\ind(1_{\ge 0}(D_a+K_a),1_{\ge 0}(D_a)) \ ,
\end{eqnarray*}
where the second equality follows from homotopy invariance and the third from the equation preceding  Prop. \ref{proj}.
Analogously $\spfl(({\tilde D}_t)_{t \in [b,b+1]})=\ind(1_{\ge 0}(D_b),1_{\ge 0}(D_b+K_b))$. 
\end{proof}

Under slightly more restricted conditions (since the previous lemma has not been proven for Hilbert $C^*$-modules -- the author did not check whether the rather complicated proof of \cite[Prop. 3.4]{le} carries over) the statement of the proposition makes sense on a Hilbert $C^*$-module and was used as a definition of the noncommutative spectral flow in \cite{wa}.

In the following we express the spectral flow in terms of a winding number.

Let $S^1=[0,1]/_{0 \sim 1}$ with the standard smooth structure.

Let ${\cal U}(H) \subset B(H)$ be the group of unitaries and let
$${\cal U}_K(H)=\{U \in {\cal U}(H) ~|~U-1 \in K(H)\} \ .$$

There is an isomorphism $$w:\pi_1({\cal U}_K(H)) \cong \bbbz$$  extending the classical winding number. In fact, if $s:S^1 \to {\cal U}_K(H)$ fulfills $(x \mapsto s(x)-1) \in C^1(S^1,l^1(H))$, where $l^1(H) \subset B(H)$ is the ideal of trace class operators endowed with the trace class norm, then $$w(s)=\frac{1}{2\pi i} \int_0^1 \Tr(s(x)^{-1} s'(x))~dx \ .$$

\begin{prop}
Let $(D_t)_{t \in [0,1]}$ be a continuous path in $\SF(H)$ with invertible endpoints. Let $\chi \in C(\bbbr)$ be a normalizing function for the map $t \mapsto D_t$ and assume that $\chi(D_0)$ and $\chi(D_1)$ are involutions. Then 
$$\spfl((D_t)_{t \in [0,1]})=w([e^{\pi i(\chi(D_t)+1)}]) \ .$$
If $(D_t)_{t \in [0,1]}$ is a continuous path in $\SF(H)$ with $D_0=D_1$ (not necessarily invertible), then this equation holds for any normalizing function of $t  \mapsto D_t$.
\end{prop}

\begin{proof} The term on the right hand side is well-defined by Lemma \ref{compcont2}.
We make use of the space $KC_{sa}(H)$, which was defined in the proof of Theorem \ref{classsp}.

Let $(D_t)_{t \in [0,1]}$ be a continuous path in $\SF(H)$ with invertible endpoints. Since for any normalizing function $\chi$ of $t \mapsto D_t$ 
$$\spfl((D_t)_{t \in [0,1]})=\spfl((\chi(D_t))_{t \in [0,1]}) \ ,$$ it is enough to prove that for any continuous path $(F_t)_{t \in [0,1]}$ in $KC_{sa}(H)$ such that $F_0, F_1$ are involutions  
$$\spfl((F_t)_{t \in [0,1]})=w([e^{\pi i(F_t +1)}]) \ . \qquad (*)$$
Both sides of this equation remain unchanged if we replace $(F_t)_{t \in [0,1]}$ by $(F_t \oplus I)_{t \in [0,1]}$, where $I$ is an involution on $H$ with infinite-dimensional eigenspaces. Hence we may assume that the eigenspaces of $F_0,F_1$ are infinite-dimensional. Then there is a unitary $U$ with $F_0=UF_1U^*$. Furthermore by the contractibility of ${\cal U}(H)$ there is a continuous path $(U_t)_{t \in [1,2]}$ of unitaries, unique up to homotopy, with $U_1=1$ and $U_2=U$. Define  $G_t=F_{2t}$ for $t \in [0,\frac 12]$ and $G_t=U_{2t}F_1U_{2t}^*$ for $t \in [\frac 12,1]$. The path $(G_t)_{t \in [0,1]}$ is a loop in $KC_{sa}(H)$ with
$$\spfl((F_t)_{t \in [0,1]})=\spfl((G_t)_{t \in [0,1]})$$
and $$w([e^{\pi i(F_t +1)}])=w([e^{\pi i(G_t +1)}]) \ .$$
Thus it is enough to prove equation $(*)$ for loops in $KC_{sa}(H)$. This will also show the second assertion of the proposition.

By homotopy invariance of the winding number and of spectral flow for loops (property (6)) and by $[S^1,KC_{sa}(H)] \cong K^1(S^1) \cong \bbbz$ it is sufficient to verify the assertion for some loop in $KC_{sa}(H)$ whose class generates $K^1(S^1)$. For example one can use the loop $(G_t)_{t\in [0,1]}$ arising as above from $F_t=-\cos(\pi t)P+(1-P)$, where $P$ is a projection whose range has dimension one. In this case equation $(*)$ is well-known since $(G_t)_{t \in [0,1]}$ is a norm-continuous path.  
\end{proof}

\textsc{Gottfried Wilhelm Leibniz Bibliothek, 
Nieders\"achsische Landesbibliothek, 
Waterloostr. 8,
30169 Hannover, Germany} 

\textsc{Email: ac.wahl@web.de}
\end{document}